\documentclass[12pt]{amsart}

\usepackage{amssymb}
\usepackage{latexsym}
\usepackage{amsmath}

\newtheorem{definition}{Definition}[section]
\newtheorem{proposition}{Proposition}[section]
\newtheorem{lemma}{Lemma}[section]
\newtheorem{theorem}{Theorem}[section]
\newtheorem{corollary}{Corollary}[section]

\begin{document}

\title{Integration over compact quantum groups}
\author{Teodor Banica}
\address{Departement of Mathematics, Universite Paul Sabatier, 118 route de Narbonne, 31062 Toulouse, France}
\email{banica@picard.ups-tlse.fr}
\thanks{Communicated by M. Kashiwara. Received November 11, 2005, Revised February 20, 2006}
\author{Beno\^\i{}t Collins $^{\dagger}$}
\address{Institut Camille Jordan,
Universit\'e Claude Bernard Lyon 1,
43 boulevard du 11 novembre 1918,
69622 Villeurbanne Cedex,
France}
\email{collins@math.univ-lyon1.fr}
\thanks{$\dagger$ Research supported by RIMS COE postdoctoral fellowship}

\begin{abstract}
We find a combinatorial formula for the Haar functional of the orthogonal and unitary quantum groups. As an application, we  consider diagonal coefficients of the fundamental representation, and we investigate their spectral measures.
\end{abstract}

\maketitle

\section*{Introduction}

A basic question in functional analysis is to find axioms for quantum groups, which ensure the existence of a Haar measure. In the compact case, this was solved by Woronowicz in the late eighties (\cite{wo1}). The Haar functional is constructed starting from an arbitrary faithful positive unital linear form $\varphi$, by taking a Cesaro limit with respect to convolution:
$$\int=\lim_{n\to\infty}\frac{1}{n}\sum_{k=1}^n\varphi^{*k}$$

The explicit computation of the Haar functional is a representation theory problem. There are basically two ideas here:

I. For a classical group the integrals can be computed by using inversion of matrices and non-crossing partitions. The idea goes back to Weingarten's work \cite{we}, and explicit formulae are found in \cite{co1}, \cite{cs}.

II. For a free quantum group the integrals of characters can be computed by using tensor categories and diagrams. The idea goes back to Woronowicz's work \cite{wo2}, and several examples are studied in \cite{ba1}, \cite{ba2}.

In this paper we find an explicit formula for the Haar functional of free quantum groups. For this purpose, we use a combination of I and II.

As an application, we consider diagonal coefficients of the fundamental representation, and we investigate their spectral measures. For instance in the orthogonal case we find a formula of type
$$\int (u_{11}+\ldots +u_{ss})^{2k}=Tr(G_{kn}^{-1}G_{ks})$$
where $G_{kn}$ is a certain Gram matrix of Temperley-Lieb diagrams. This enables us to find several partial results regarding the law of $u_{11}$.

The interest here is that knowledge of the law of $u_{11}$ would be the first step towards finding a model for the orthogonal quantum group. That is, searching for an explicit operator $U_{11}$ doing what the abstract operator $u_{11}$ does would be much easier once we know its law.

As a conclusion, we can state some precise problems. In the orthogonal case the question is to find the real measure $\mu$ satisfying 
$$\int x^{2k}\,d\mu(x)=Tr(G_{kn}^{-1}G_{ks})$$
and we have a similar statement in the unitary case.

An answer to these questions would no doubt bring new information about free quantum groups. But this requires a good knowledge of combinatorics of Gram matrices, that we don't have so far.

The whole thing is probably related to questions considered by Di Francesco, Golinelli and Guitter, in connection with the meander problem. In \cite{df}, \cite{df+} they find a formula for the determinant of $G_{kn}$, but we don't know yet how to apply their techniques to our situation.

Finally, let us mention that techniques in this paper apply as well to the quantum symmetric group and its versions, whose corresponding Hom spaces are known to be described by Temperley-Lieb diagrams (\cite{ba2}, \cite{wa2}). 
This will be discussed in a series of papers, the first of which is in preparation (\cite{bc}).

The paper is organised as follows. 1, 2, 3 are preliminary sections on the orthogonal quantum group. In 4, 5, 6, 7, 8 we establish the orthogonal integration formula, then we apply it to diagonal coefficients, and then to coefficients of type $u_{11}$, with a separate discussion of the case $n=2$. In 9 we find similar results for the unitary quantum group.

\subsection*{Acknowledgements}  
We would like to express our deepest gratitude to the referee, for a careful reading of the manuscript.

\section{The orthogonal quantum group}

In this section we present a few basic facts regarding the universal algebra $A_o(n)$. This algebra appears in Wang's thesis (see \cite{wa1}).

For a square matrix $u=u_{ij}$ having coefficients in a ${\mathcal C}^*$-algebra, we use the notations $\bar{u}=u_{ij}^*$, $u^t=u_{ji}$ and $u^*=u_{ji}^*$.

A matrix $u$ is called orthogonal if $u=\bar{u}$ and $u^t=u^{-1}$.

\begin{definition}
$A_o(n)$ is the ${\mathcal C}^*$-algebra generated by $n^2$ elements $u_{ij}$, with relations making $u=u_{ij}$ an orthogonal matrix.
\end{definition}

In other words, we have the following universal property. For any pair $(B,v)$ consisting of a ${\mathcal C}^*$-algebra $B$ and an orthogonal matrix $v\in M_n(B)$, there is a unique morphism of ${\mathcal C}^*$-algebras
$$A_o(n)\to B$$
mapping $u_{ij}\to v_{ij}$ for any $i,j$. The existence and uniqueness of such a universal pair $(A_o(n),u)$ follow from standard ${\mathcal C}^*$-algebra results.

\begin{proposition}
$A_o(n)$ is a Hopf ${\mathcal C}^*$-algebra, with comultiplication, counit and antipode given by the formulae
$$\Delta(u_{ij})=\sum_{k=1}^n u_{ik}\otimes u_{kj}$$
$$\varepsilon(u_{ij})=\delta_{ij}$$
$$S(u_{ij})=u_{ji}$$
which express the fact that $u$ is a $n$-dimensional corepresentation.
\end{proposition}

These maps are constructed by using the universal property of $A_o(n)$, and verification of Woronowicz's axioms in \cite{wo1} is straightforward. As an example, the counit $\varepsilon :A_o(n)\to {\mathbb C}$ is constructed by using the fact that $1_n=\delta_{ij}$ is an orthogonal matrix over the algebra ${\mathbb C}$.

Observe that the square of the antipode is the identity:
$$S^2=id$$

The motivating fact about $A_o(n)$ is a certain analogy with ${\mathbb C}(O(n))$. The coefficients $v_{ij}$ of the fundamental representation of $O(n)$ form an orthogonal matrix, and we have the following presentation result.

\begin{proposition}
${\mathcal C}(O(n))$ is the commutative ${\mathcal C}^*$-algebra generated by $n^2$ elements $v_{ij}$, with relations making $v=v_{ij}$ an orthogonal matrix.
\end{proposition}

Observe in particular that we have a morphism of ${\mathcal C}^*$-algebras
$$A_o(n)\to {\mathcal C}(O(n))$$
mapping $u_{ij}\to v_{ij}$ for any $i,j$. The above formulae of $\Delta,\varepsilon,S$ show that this is a Hopf algebra morphism. We get an isomorphism
$$A_o(n)/I= {\mathcal C}(O(n))$$
where $I$ is the following ideal:
$$I=<[u_{ij},u_{kl}]=0\mid i,j,k,l>$$

This is usually called commutator ideal, because the quotient by it is the biggest commutative quotient.

This result is actually not very relevant, because $A_o(n)$ has many other quotients. Consider for instance the group ${\mathbb Z}_2=\{ 1,g\}$. The equality $g=g^{-1}$ translates into the equality
$$g=g^*=g^{-1}$$
at the level of the group algebra ${\mathbb C}^*({\mathbb Z}_2)$, which tells us that the $1\times 1$ matrix $g$ is orthogonal.

Now by taking $n$ free copies of ${\mathbb Z}_2$, we get the following result.

\begin{proposition}
${\mathcal C}^*({\mathbb Z}_2^{*n})$ is the ${\mathcal C}^*$-algebra generated by $n$ elements $g_i$, with relations making $g=diag(g_1,\ldots ,g_n)$ an orthogonal matrix.
\end{proposition}

In particular we have a morphism of ${\mathcal C}^*$-algebras
$$A_o(n)\to {\mathcal C}^*({\mathbb Z}_2^{*n})$$
mapping $u_{ij}\to g_{ij}$ for any $i,j$. The above formulae of $\Delta,\varepsilon,S$ show that this is a Hopf algebra morphism. We get an isomorphism
$$A_o(n)/J= {\mathcal C}^*({\mathbb Z}_2^{*n})$$
where $J$ is the following ideal:
$$J=<u_{ij}=0\mid i\neq j>$$

This can be probably called cocommutator ideal, because the quotient by it is the biggest cocommutative quotient.

As a conclusion here, best is to draw a diagram.

\begin{theorem}
We have surjective morphisms of Hopf ${\mathcal C}^*$-algebras
$$\begin{matrix}
&&A_o(n)&&\cr
&\swarrow&&\searrow&\cr
{\mathcal C}(O(n))&&&&{\mathcal C}^*({\mathbb Z}_2^{*n})
\end{matrix}$$
obtained from the universal property of $A_o(n)$.
\end{theorem}

This diagram is to remind us that $A_o(n)$ is at the same time a non-commutative version of ${\mathcal C}(O(n))$, and a non-cocommutative version of ${\mathcal C}^*({\mathbb Z}_2^{*n})$. We say that it is a free version of both.

\section{Analogy with $SU(2)$}

The study of $A_o(n)$ is based on a certain similarity with ${\mathcal C}(SU(2))$. The fundamental corepresentation of ${\mathcal C}(SU(2))$ is given by
$$w=\begin{pmatrix}a&b\cr -\bar{b}&\bar{a}\end{pmatrix}$$
with $|a|^2+|b|^2=1$. This is of course a unitary matrix, which is not orthogonal. However, $w$ and $\bar{w}$ are related by the formula 
$$\begin{pmatrix}a&b\cr -\bar{b}&\bar{a}\end{pmatrix}\begin{pmatrix}0&1\cr -1&0\end{pmatrix}=\begin{pmatrix}0&1\cr -1&0\end{pmatrix}\begin{pmatrix}\bar{a}&\bar{b}\cr -b&a\end{pmatrix}$$
which is a twisted self-conjugation condition of type
$$w=r\bar{w}r^{-1}$$
where $r$ is the following matrix:
$$r=\begin{pmatrix}0&1\cr -1&0\end{pmatrix}$$

One can show that unitarity plus this condition are in fact the only ones, in the sense that we have the following presentation result.

\begin{proposition}
${\mathcal C}(SU(2))$ is the ${\mathcal C}^*$-algebra generated by $4$ elements $w_{ij}$, with the relations $w=r\bar{w}r^{-1}=\mbox{unitary}$, where $w=w_{ij}$.
\end{proposition}

This is to be compared with the definition of $A_o(n)$, which can be written in the following way.
$$A_o(n)={\mathcal C}^*\left\{ (u_{ij})_{ij=1,\ldots ,n}\,\,|\,\,u=\bar{u}=\mbox{unitary}\right\}$$

We see that what makes the difference between the two matrices $v_1=u$ and $v_2=w$ is possibly their size, plus the value of a scalar matrix $r$ intertwining $v$ and $\bar{v}$.

This leads to the conclusion that $A_o(n)$ should be a kind of deformation of ${\mathcal C}(SU(2))$. Here is a precise result in this sense.

\begin{theorem}
We have an isomorphism
$$A_o(2)={\mathcal C}(SU(2))_{-1}$$
where the algebra on the right is the specialisation at $\mu=-1$ of the algebra ${\mathcal C}(SU(2))_\mu$ constructed by Woronowicz in  \cite{wo0}.
\end{theorem}

This result, pointed out in \cite{ba1}, is clear from definitions.

We should mention here that the parameter $\mu\in {\mathbb R}-\{ 0\}$ used by Woronowicz in \cite{wo0} is not a particular case of the parameter $q\in {\mathbb C}-\{ 0\}$ used in the quantum group literature. In fact, we have the formula
$$\mu =\tau q^2$$
where $q>0$ is the usual deformation parameter, and where $\tau=\pm 1$ is the twist, constructed by Kazhdan and Wenzl in \cite{kw}. In particular the value $\mu =-1$ corresponds to the values $q=1$ and $\tau =-1$.

Finally, let us mention that theorem 2.1 follows via a change of variables from the general formula
$$A_o\begin{pmatrix}0&1\cr -\mu^{-1}&0\end{pmatrix}={\mathcal C}(SU(2))_\mu$$
where the algebra on the left is constructed in the following way:
$$A_o(r)={\mathcal C}^*\left\{ (u_{ij})_{ij=1,\ldots ,n}\,\,|\,\,u=r\bar{u}r^{-1}=\mbox{unitary}\right\}$$

See \cite{bdv} for more on parametrisation of algebras of type $A_o(r)$.

\section{Diagrams}

The main feature of the fundamental representation of $SU(2)$ is that commutants of its tensor powers are Temperley-Lieb algebras:
$$End(w^{\otimes k})=TL(k)$$

This equality is known to hold in fact for the fundamental corepresentation of any ${\mathcal C}(SU(2))_\mu$, as shown by Woronowicz in \cite{wo0}.

The same happens for $A_o(n)$, as pointed out in \cite{ba1}. We present now a proof of this fact, a bit more enlightening than the original one. For yet another proof, see Yamagami (\cite{ya1}, \cite{ya2}).

\begin{definition}
The set of Temperley-Lieb diagrams $D(k,l)$ consists of diagrams formed by an upper row of $k$ points, a lower row of $l$ points, and of $(k+l)/2$ non-crossing strings joining pairs of points.
\end{definition}

In this definition, for $k+l$ odd we have $D(k,l)=\emptyset$. Also, diagrams are taken of course up to planar isotopy.

It is convenient to summarize this definition as
$$D(k,l)=\left\{ \begin{matrix}\cdot\,\cdot\,\cdot & \leftarrow &
    k\,\,\mbox{points}\cr W & \leftarrow &
    (k+l)/2\mbox{ strings}\cr \cdot\,\cdot\,\cdot\,\cdot\,\cdot& \leftarrow &
    l\,\,\mbox{points}\end{matrix}\right\}$$
where capital letters denote diagrams formed by non-crossing strings.

\begin{definition}
The operation on diagrams given by
$$\begin{matrix}\cdot\,\cdot\,\cdot\cr
W\cr ||\cr A\cr
\cdot\,\cdot\,\cdot\,\cdot\,\cdot\end{matrix}
\ \ \rightarrow\ \ 
\begin{matrix}
\Cap\cr
A\, M\cr
\cdot\,\cdot\,\cdot\,\cdot\,\cdot\,\cdot\,\cdot\,\cdot\end{matrix}$$
is an identification $D(k,l)\simeq D(0,k+l)$, called Frobenius isomorphism.
\end{definition}

Observe in particular the identification at $k=l$, namely
$$D(k,k)\simeq D(0,2k)$$
where at left we have usual Temperley-Lieb diagrams, 
$$D(k)=\left\{ \begin{matrix}\cdot\,\cdot\,\cdot & \leftarrow &
    k\,\,\mbox{points}\cr W & \leftarrow &
    k\mbox{ strings}\cr \cdot\,\cdot\,\cdot\,\cdot\,\cdot& \leftarrow &
    k\,\,\mbox{points}\end{matrix}\right\}$$
and at right we have non-crossing partitions of $1,\ldots ,2k$:
$$NC(2k)=\left\{ \begin{matrix}
    W & \leftarrow &
    k\mbox{ strings}\cr \cdot\,\cdot\,\cdot\,\cdot\,\cdot& \leftarrow &
    2k\,\,\mbox{points}\end{matrix}\right\}$$

It is convenient to reformulate the above Frobenius isomorphism by using these notations, and to use it as an equality.

\begin{definition}
We use the Frobenius identification 
$$D(k)=NC(2k)$$
between usual Temperley-Lieb diagrams and non-crossing partitions.
\end{definition}

Consider now the vector space where $v$ acts, namely
$$V={\mathbb C}^n$$
and denote by $e_1,\ldots ,e_n$ its standard basis. Each diagram $p\in D(k,l)$ acts on tensors according to the formula
$$p(e_{i_1}\otimes\ldots\otimes e_{i_k})=\sum_{j_1\ldots j_l}\begin{pmatrix}i_1\ldots i_k\cr p\cr j_1\ldots j_l\end{pmatrix}e_{j_1}\otimes\ldots\otimes e_{j_l}$$
where the middle symbol is $1$ if all strings of $p$ join pairs of equal indices, and is $0$ if not. Linear maps corresponding to different diagrams can be shown to be linearly independent provided that $n\geq 2$, and this gives an embedding $$TL(k,l)\subset Hom(V^{\otimes k},V^{\otimes l})$$
where $TL(k,l)$ is the abstract vector space spanned by $D(k,l)$. This is easy to check by using positivity of the trace, see for instance \cite{ba2}.

\begin{theorem}
We have an equality of vector spaces
$$Hom(u^{\otimes k},u^{\otimes l})=TL(k,l)$$
where 
$Hom(u^{\otimes k},u^{\otimes l})$ is the subalgebra of 
$Hom(V^{\otimes k},V^{\otimes l})$ of $A_o(n)$-equivariant endomorphisms.
\end{theorem}

\begin{proof}
We use tensor categories with suitable positivity properties, as axiomatized by Woronowicz in \cite{wo2}.

The starting remark is that for a unitary matrix $u$, the fact that $u$ is orthogonal is equivalent to the fact that the vector
$$\xi=\sum_k e_k\otimes e_k$$
is fixed by $u^{\otimes 2}$, in the sense that we have the following equality:
$$u^{\otimes 2}(\xi\otimes 1)=\xi\otimes 1$$

This follows by writing down relations for both conditions on $u$. Now in terms of the linear map $E:{\mathbb C}\to V^{\otimes 2}$ given by
$$E(1)=\xi$$
we have the following equivalent condition:
$$E\in Hom(1,u^{\otimes 2})$$

On the other hand, $E$ is nothing but the operator corresponding to the semicircle in $D(0,2)$:
$$E=\cap$$

Summing up, $A_o(n)$ is the universal ${\mathcal C}^*$-algebra generated by entries of a unitary $n\times n$ matrix $u$, satisfying the following condition:
$$\cap\in Hom(1,u^{\otimes 2})$$

In terms of tensor categories, this gives the equality
$$<\cap>=\{ Hom(u^{\otimes k},u^{\otimes l})\mid k,l\}$$
where the category on the left is the one generated by $\cap$, meaning the smallest one satisfying Woronowicz's axioms in \cite{wo2}, and containing $\cap$.

Woronowicz's operations are the composition, tensor product and conjugation. At level of Temperley-Lieb diagrams, these are easily seen to correspond to horizontal concatenation, vertical concatenation and upside-down turning of diagrams. Since all Temperley-Lieb diagrams can be obtained from $\cap$ via these operations, we get the equality
$$<\cap>=\{ TL(k,l)\mid k,l\}$$
which together with the above equality gives the result.
\end{proof}

Observe that the ingredients of this proof are Woronowicz's Tannakian duality, plus basic facts concerning Temperley-Lieb diagrams. For a more detailed application of Tannakian duality, in a similar situation, see \cite{ba2}. As for Temperley-Lieb diagrams, what we use here is the tensor planar algebra, constructed by Jones in \cite{j}.

\section{Integration formula}

In this section we find a formula for the Haar functional of $A_o(n)$. This is a certain linear form, denoted here as an integral
$$\int :A_o(n)\to {\mathbb C}$$
and whose fundamental property is the following one.

\begin{definition}
The Haar functional of $A_o(n)$ is the positive linear unital form satisfying the bi-invariance condition
$$\left( id\otimes \int\right)\Delta(a) =\left(\int\otimes id\right)\Delta(a) =\int a$$
whose existence and uniqueness is shown by Woronowicz in \cite{wo1}.
\end{definition}

For the purposes of this paper, we just need the following property: for a unitary corepresentation $r\in End(H)\otimes A_o(n)$, the operator
$$P=\left( id\otimes\int\right)r$$
is the orthogonal projection onto the space of fixed points of $r$. This space is in turn defined as
$$Hom(1,r)=\{ x\in H \mid r(x)=x\otimes 1\}$$
and the whole assertion is proved in \cite{wo1}.

The integration formula involves scalar matrices $G_{kn}$ and $W_{kn}$, introduced in the following way.

\begin{definition}
The Gram and Weingarten matrices are given by
$$G_{kn}(p,q)=n^{l(p,q)}$$
$$W_{kn}=G_{kn}^{-1}$$
where $l(p,q)$ is the number of loops obtained by closing the composed diagram $p^*q$ for
$p,q\in D(k)$.
\end{definition}

The fact that $G_{kn}$ is indeed a Gram matrix comes from the equality
$$G_{kn}(p,q)=<p,q>$$
where $p,q$ are regarded as operators on the Hilbert space $V^{\otimes k}$, with $V={\mathbb C}^n$, and where the scalar product on $V$ is the usual one.
Alternatively, $<p,q>$ can be understood as the value of the Markov trace of $p^*q$
in the Temperley-Lieb algebra. %

As for $W_{kn}$, we will see that this is a quantum analogue of the matrix constructed by Weingarten in \cite{we}.

For a diagram $p\in D(k)$ and a multi-index $i=(i_1\ldots i_{2k})$ we use the notation
$$\delta_{pi}=\begin{pmatrix}i_{2k}\ldots i_{k+1}\cr p\cr i_1\ldots i_{k}\end{pmatrix}$$
where, as usual, the symbol on the right is $1$ if all strings of $p$ join pairs of equal indices, and $0$ if not. This is the same as the notation
$$\delta_{pi}=\begin{pmatrix}p\cr i_1\ldots i_{2k}\end{pmatrix}$$
where $p$ is regarded now as a non-crossing partition, via the Frobenius identification in definition 3.3.

\begin{theorem}
The Haar functional of $A_o(n)$ is given by
$$\int u_{i_{1}j_{1}}\ldots u_{i_{2k}j_{2k}}=\sum_{pq}\delta_{pi}\delta_{qj}W_{kn}(p,q)$$
$$\int u_{i_{1}j_{1}}\ldots u_{i_{2k+1}j_{2k+1}}=0$$
where the sum is over all pairs of diagrams $p,q\in D(k)$.
\end{theorem}

\begin{proof}
We have to compute the linear maps
$$E(e_{i_{1}}\otimes \ldots \otimes e_{i_{l}})=
\sum_{j_{1}\ldots j_{l}}
e_{j_{1}}\otimes \ldots \otimes e_{j_{l}}
\int u_{i_{1}j_{1}}\ldots u_{i_{l}j_{l}}$$
which encode all integrals in the statement.

In case $l=2k$ is even we use the fact that $E$ is the orthogonal projection onto $End(u^{\otimes k})$. With the notation
$$\Phi (x)=\sum_p<x,p>p$$
we have $E=W\Phi$, where $W$ is the inverse on $TL(k)$ of the restriction of $\Phi$. But this restriction is the linear map given by $G_{kn}$, so $W$ is the linear map given by $W_{kn}$. This gives the first formula. 

In case $l$ is odd we use the automorphism $u_{ij}\to -u_{ij}$ of $A_o(n)$. From $E=(-1)^lE$ we get $E=0$, which proves the second formula.
\end{proof}

\section{Diagonal coefficients}

The law of a self-adjoint element $a\in A_o(n)$ is the real probability measure $\mu$ given by
$$\int\varphi(x)\,d\mu(x)=\int \varphi (a)$$
for any continuous function $\varphi:{\mathbb R}\to{\mathbb C}$. 
As for any bounded probability measure, $\mu$ is uniquely determined by its moments. These are the numbers
$$\int x^k\,d\mu(x)=\int a^k$$
with $k=1,2,3,\ldots$, also called moments of $a$.

We are particularly interested in the following choice of $a$.

\begin{definition}
The $o_{sn}$ variable is given by
$$o_{sn}=u_{11}+\ldots +u_{ss}$$
where $u$ is the fundamental corepresentation of $A_o(n)$.
\end{definition}

The motivating fact here is that all coefficients $u_{ii}$ have the same law. This is easily seen by using automorphisms of $A_o(n)$ of type
$$\sigma:u\to pup^{-1}$$
where $p$ is a permutation matrix. This common law, whose knowledge might be the first step towards finding an explicit model for $A_o(n)$, is the law of $o_{1n}$.

The idea of regarding $o_{1n}$ as a specialisation of $o_{sn}$ comes from the fact that $o_{nn}$ is a well-known variable, namely the semicircular one. This is known from \cite{ba1}, and is deduced here from the following result.

\begin{theorem}
The even moments of the $o_{sn}$ variable are given by
$$\int o_{sn}^{2k}=Tr(G_{kn}^{-1}G_{ks})$$
and the odd moments are all equal to $0$.
\end{theorem}

\begin{proof}
The first assertion follows from theorem 4.1,
\begin{eqnarray*}
\int o_{sn}^{2k}
&=&\int (u_{11}+\ldots +u_{ss})^{2k}\cr
&=&\sum_{a_1=1}^{s}\ldots\sum_{a_{2k}=1}^s\int u_{a_1a_1}\ldots u_{a_{2k}a_{2k}}\cr
&=&\sum_{a_1=1}^{s}\ldots\sum_{a_{2k}=1}^s\sum_{p,q\in D(k)}\delta_{pa}\delta_{qa}W_{kn}(p,q)\cr
&=&\sum_{p,q\in D(k)}W_{kn}(p,q)\sum_{a_1=1}^{s}\ldots\sum_{a_{2k}=1}^s\delta_{pa}\delta_{qa}\cr
&=&\sum_{p,q\in D(k)}W_{kn}(p,q)G_{ks}(q,p)\cr
&=&Tr(W_{kn}G_{ks})
\end{eqnarray*}
and from the equality $W_{kn}=G_{kn}^{-1}$. As for the assertion about odd moments, this follows as well from theorem 4.1.
\end{proof}

As a first application, we get another proof for the fact that $o_{nn}$ is semicircular. The semicircle law has density
$$d\mu (x)=\frac{1}{2\pi}\sqrt{4-x^2}\,dx$$
on $[-2,2]$, and $0$ elsewhere. A variable having this law is called semicircular. The even moments of $\mu$ are the Catalan numbers
$$C_k=\frac{1}{k+1}\begin{pmatrix}2k\cr k\end{pmatrix}$$
and the odd moments are all equal to $0$. See \cite{vdn}.

\begin{corollary}
The variable $o_{nn}$ is semicircular.
\end{corollary}

\begin{proof}
The even moments of $o_{nn}$ are the Catalan numbers
\begin{eqnarray*}
\int o_{nn}^{2k}
&=&Tr(G_{kn}^{-1}G_{kn})\cr
&=&Tr(1)\cr
&=&\# D(k)\cr
&=&C_k
\end{eqnarray*}
hence are equal to the even moments of the semicircle law. As for odd moments, they are $0$ for both $o_{nn}$ and for the semicircle law.
\end{proof}

The second application brings some new information about $A_o(n)$.

\begin{corollary}
The variable $(n/s)^{1/2}\,o_{sn}$ is asymptotically semicircular as $n\to\infty$.
\end{corollary}

\begin{proof}
We have $G_{kn}(p,q)=n^k$ for $p=q$, and $G_{kn}(p,q)\leq n^{k-1}$ for $p\neq q$. Thus with $n\to\infty$ we have $G_{kn}\sim n^k1$, which gives
\begin{eqnarray*}
\int o_{sn}^{2k}
&=&Tr(G_{kn}^{-1}G_{ks})\cr
&\sim&Tr((n^k1)^{-1} G_{ks})\cr
&=&n^{-k}Tr(G_{ks})\cr
&=&n^{-k}s^k\# D(k)\cr
&=&n^{-k}s^kC_k
\end{eqnarray*}
which gives the convergence in the statement, for even moments. As for odd ones, they are all $0$, so we have convergence here as well.
\end{proof}

\section{Asymptotic freeness}

We know from corollary 5.2 that the variable $n^{1/2}o_{1n}$ is asymptotically semicircular. Together with the observation after definition 5.1, this shows that the normalised generators
$$\{n^{1/2}u_{ij}\}_{i,j=1,\ldots ,n}$$
of $A_o(n)$ become asymptotically semicircular as $n\to\infty$. Here we assume that $i,j$ are fixed, say $i,j\leq s$ and the limit is over $n\geq s$.

This result might be useful when looking for explicit models for $A_o(n)$. Here is a more precise statement in this sense.

\begin{theorem}
The elements $(n^{1/2}u_{ij})_{i,j=1,\ldots ,s}$ of $A_o(n)$ with $n\geq s$ become asymptotically free and semicircular as $n\to\infty$.
\end{theorem}

\begin{proof}
The joint moments of a free family of semicircular elements are computed by using the fact that the second order free cumulant is one, and the other ones are zero.
Therefore for a free family of semicircular variables 
$x_{1},\ldots ,x_{k}$, an integral of type
$$\int x_{i_{1}}\ldots x_{i_{l}}$$ 
is zero if $l$ is odd, and is the sum of matching non-crossing pair partitions if $l$ is
even. This is a free version of Wick theorem; see Speicher (\cite{s1}) for details. Now when computing 
$$n^{k}\int u_{i_{1}j_{1}}\ldots u_{i_{2k}j_{2k}}$$
by using Theorem 4.1, observe that 
$$n^{k}W_{kn}(p,p)\to 1$$
$$n^{k}W_{kn}(p,q)\to 0$$
as $n\to \infty$, whenever $p\neq q$. This completes the proof.
\end{proof}

Observe that theorem 6.1 is indeed stronger than corollary 5.2: it is known that, with suitable normalisations, a sum of free semicircular variables is semicircular. See \cite{vdn}.

\section{Second order results}

A basic problem regarding the algebra $A_o(n)$ is to find the law of coefficients $u_{ij}$. This is the law of the variable $o_{1n}$, as defined in previous section, with moments given by
$$\int o_{1n}^{2k}=\sum_{p,q\in D(k)} W_{kn}(p,q).$$

We know from corollary 5.2 that, under a suitable normalisation, these moments converge with $n\to\infty$ to those of the semicircle law. In this section we find a power series expansion of $W_{kn}$, which can be used for finding higher order results about the law of $o_{1n}$.

Observe first that the integer-valued function
$$d(p,q)=k-l(p,q)$$
is a distance on the space $D(k)$. Indeed, it can be shown by induction
that if $p\neq q$, $d(p,q)$ is the minimal number $l$ such that there exists 
$p_1,\ldots ,p_l$ satisfying $p_1=p$, $p_l=q$, and for each 
pair $\{p_i,p_{i+1}\}$, $p_i,p_{i+1}$ have all strings identical except two of them.
We call this distance ``loop distance".

\begin{proposition}
The Gram matrix is given by
$$n^{-k}G_{kn}(p,q)=n^{-d(p,q)}$$
where $d$ is the loop distance on $D(k)$.
\end{proposition}

\begin{proof}
This is clear from definitions.
\end{proof}

In other words, the matrix $n^{-k}G_{kn}$ is an entry-wise exponential of the distance matrix of $D(k)$. This exponential can be inverted by using paths on $D(k)$. Such a path is a sequence of elements of the form:
$$p_{0}\neq p_{1}\neq \ldots \neq p_{l-1}
\neq p_{l}$$

We call this sequence path from $p_0$ to $p_l$.

\begin{definition}
The distance along a path $P=p_0,\ldots ,p_l$ is the number
$$d(P)=d(p_0,p_1)+\ldots +d(p_{l-1},p_l)$$
and the length of such a path is the number $l(P)=l$.
\end{definition}

Observe that a length $0$ path is just a point, and the distance along such a path is $0$.

With these definitions, we have a power series expansion in $n^{-1}$ for the Weingarten matrix.

\begin{proposition}
The Weingarten matrix is given by
$$n^kW_{kn}(p,q)=\sum_P(-1)^{l(P)}n^{-d(P)}$$
where the sum is over all paths from $p$ to $q$.
\end{proposition}

\begin{proof}
For $n$ large enough we have the following computation.
\begin{eqnarray*}
n^{k}W_{kn}
&=&(n^{-k}G_{kn})^{-1}\cr
&=&(1-(1-n^{-k}G_{kn}))^{-1}\cr
&=&1+\sum_{l=1}^\infty (1-n^{-k}G_{kn})^l
\end{eqnarray*}

We know that $G_{kn}$ has $n^k$ on its diagonal, so $1-n^{-k}G_{kn}$ has $0$ on the diagonal, and its $l$-th power is given by
\begin{eqnarray*}
(1-n^{-k}G_{kn})^{l}(p,q)
&=&\sum_P\prod_{i=1}^l(1-n^{-k}G_{kn})(p_{i-1},p_i)\cr
&=&\sum_P\prod_{i=1}^l-n^{-d(p_{i-1},p_i)}\cr
&=&(-1)^{l}\sum_Pn^{-d(P)}
\end{eqnarray*}
with $P=p_0,\ldots ,p_l$ ranging over all length $l$ paths from $p_0=p$ to $p_l=q$. Together with the first formula, this gives
$$n^kW_{kn}(p,q)=\delta_{pq}+\sum_P(-1)^{l(P)}n^{-d(P)}$$
where the sum is over all paths between $p$ and $q$, having length $l\geq 1$. But the leading term $\delta_{pq}$ can be added to the sum, by enlarging it to length $0$ paths, and we get the formula in the statement.
\end{proof}

In terms of moments of $o_{1n}$, we get the following power series expansion in $n^{-1}$.

\begin{proposition}
The moments of $n^{1/2}o_{1n}$ are given by
$$\int \left( n^{1/2}o_{1n}\right)^{2k}=\sum_{d=0}^\infty (E^k_d-O^k_d)n^{-d}$$
where $E^k_d,O^k_d$ count even and odd length paths of $D(k)$ of distance $d$.
\end{proposition}

\begin{proof}
>From theorem 5.1 and proposition 7.2 we get
\begin{eqnarray*}
n^k\int o_{1n}^{2k}
&=&\sum_P (-1)^{l(P)}n^{-d(P)}
\end{eqnarray*}
where the sum is over all paths in $D(k)$. This is a series in $n^{-1}$, whose $d$-th coefficient is the sum of numbers $(-1)^{l(P)}$, given by $E^k_d-O^k_d$.
\end{proof}

We have now all ingredients for computing the second order term of the law of $n^{1/2}o_{1n}$. Consider the formula
$$\int\frac{1}{1-z(n^{1/2}o_{1n})}=\sum_{k=0}^\infty z^k \int (n^{1/2}o_{1n})^{k}$$
valid for $z$ small complex number, or for $z$ formal variable. The left term is the Stieltjes transform of the law of $n^{1/2}o_{1n}$, and we have the following power series expansion of it, when $z$ is formal.

\begin{theorem}
We have the formal estimate
\begin{eqnarray*}
\sum_{k=0}^\infty z^k \int (n^{1/2}o_{1n})^{k}
&=&\frac{2}{1+\sqrt{1-4z^2}}\cr
&+&n^{-1}
\frac{32 z^4}{(1+\sqrt{1-4z^2})^4\sqrt{1-4z^2}}
\cr
&+&O(n^{-2})
\end{eqnarray*}
where $O(n^{-2})$ should be understood coefficient-wise. 
\end{theorem}

\begin{proof}
We use proposition 7.3. Since paths of distance $0$ are of length $0$ and correspond to points of $D(k)$, the leading terms of the series of moments of $n^{1/2}o_{1n}$ are the Catalan numbers
$$E^k_0-O^k_0=E_0^k=\# D(k)=C_k$$
which are the moments of the semicircle law.

The next terms come from paths of distance $1$. Such a path must be of the form $P=p,q$ with $d(p,q)=1$, and has length $1$. It follows that the second terms we are interested in are given by
$$E_1^k-O_1^k=-O_1^k=-N_k$$
where $N_k$ counts neighbors in $D(k)$, meaning pairs of diagrams $(p,q)$ at distance $1$. This situation happens when all blocks of $p$ and $q$ are the same, except for two blocks of $p$ and two blocks of $q$, which do not match with corresponding blocks of $q$ and $p$. In such a situation, these four blocks yield a circle. 

Consider the generating series of numbers $C_k$ and $N_k$:
$$C(z)=\sum_{k=0}^\infty C_kz^{2k}$$
$$N(z)=\sum_{k=0}^\infty N_kz^{2k}$$

In order to make an effective enumeration of $N_k$ using power series tools, we need to make some observations:

1. The circle given by non-matching blocks of $p$ and $q$ intersects in four points
the set of $2k$ points on which elements of $D(k)$ are drawn. For each choice of four such points there are two possible circles, explaining the $2$ factor appearing
in the functional equation below.

2. There is a symmetry by circular permutation in the enumeration problem of $N_k$.
 
These observations give the following equation:
$$N(z)=2z^4C(z)^3\left(C(z)+z\frac{d}{dz}\,C(z)\right)$$

On the other hand, the generating series of Catalan numbers is
$$C(z)=\frac{2}{1+\sqrt{1-4z^2}}$$
where the square root is defined as analytic continuation on $\mathbb{C}-\mathbb{R}_{-}$ of the positive function $t\to \sqrt{t}$ on $\mathbb{R}_+^*$. We get
$$N(z)=\frac{32 z^4}{(1+\sqrt{1-4z^2})^4\sqrt{1-4z^2}}$$
which completes the proof.
\end{proof}

%
%
%
%
%
%
%

\section{The case $n=2$}

We end the study of $o_{1n}$ with a complete computation for $n=2$. The formula in this section is probably known to specialists, because $A_o(2)$ is one of the much studied deformations of ${\mathbb C}(SU(2))$, but we were unable to find the right bibliographical reference for it.

\begin{lemma}
We have the equalities
$$u_{11}^{2}+u_{12}^{2}=1$$
$$[u_{12},u_{11}^{2}]=0$$
where $v$ is the fundamental corepresentation of $A_o(2)$.
\end{lemma}

\begin{proof}
The first equality comes from the fact that $u$ is orthogonal. The second one comes from the computation
\begin{eqnarray*}
u_{12}u_{11}^{2}-u_{11}^{2}u_{12}
&=&u_{12}(1-u_{12}^{2})-(1-u_{12}^{2})u_{12}\cr
&=&u_{12}-u_{12}^{3}-u_{12}+u_{12}^3\cr
&=&0
\end{eqnarray*}
where we use twice the first equality.
\end{proof}

\begin{theorem}
For the generators $u_{ij}$ of the algebra $A_o(2)$, the law of each $u_{ij}^2$ is the uniform measure on $[0,1]$.
\end{theorem}

\begin{proof}
As explained after definition 5.1, we may assume $i=j=1$. Let $D=D(k)$. We use the partition
$$D=D_1\sqcup\ldots\sqcup D_k$$
where $D_i$ is the set of of diagrams such that a string joins $1$ with $2i$.

By applying twice theorem 4.1, then by using several times lemma 8.1, we have the following computation.
\begin{eqnarray*}
\int u_{11}^{2k}
&=&\sum_{p,q\in D}W_{k2}(p,q)\cr
&=&\sum_{l=1}^{k}\sum_{p\in D}\sum_{q\in D_l}W_{k2}(p,q)\cr
&=&\sum_{l=1}^{k}\int u_{12}u_{11}^{2l-2}u_{12}u_{11}^{2k-2l}\cr
&=&\sum_{l=1}^{k}\int u_{12}^2u_{11}^{2k-2}\cr
&=&\sum_{l=1}^{k}\int (1-u_{11}^2)u_{11}^{2k-2}\cr
&=&k\int u_{11}^{2k-2}-k\int u_{11}^{2k}
\end{eqnarray*}

Rearranging terms gives the formula
$$(k+1)\int u_{11}^{2k}=k\int u_{11}^{2k-2}$$
and we get by induction on $k$ the value of all moments of $u_{11}^2$:
$$\int u_{11}^{2k}=\frac{1}{k+1}$$

But these numbers are known to be the moments of the uniform measure on $[0,1]$, and we are done.
\end{proof}

\section{The unitary quantum group}

In this section we study the Haar functional of the universal algebra $A_u(n)$. This algebra appears in Wang's thesis (see \cite{wa1}).

\begin{definition}
$A_u(n)$ is the ${\mathcal C}^*$-algebra generated by $n^2$ elements $v_{ij}$, with relations making $v=v_{ij}$ and $v^t=v_{ji}$ unitary matrices.
\end{definition}

It follows from definitions that $A_u(n)$ is a Hopf ${\mathcal C}^*$-algebra. The comultiplication, counit and antipode are given by the formulae
$$\Delta(v_{ij})=\sum_{i=1}^n v_{ik}\otimes v_{kj}$$
$$\varepsilon(v_{ij})=\delta_{ij}$$
$$S(v_{ij})=v_{ji}^*$$
which express the fact that $v$ is an $n$-dimensional corepresentation.

The motivating fact about $A_u(n)$ is an analogue of theorem 1.1, involving the unitary group $U(n)$ and the free group $F_n$.
$$\begin{matrix}
&&A_u(n)&&\cr
&\swarrow&&\searrow&\cr
{\mathbb C}(U(n))&&&&{\mathbb C}^*(F_n)
\end{matrix}$$

We already know that this kind of result might not be very relevant. This is indeed the case, so we switch to computation of commutants. For this purpose, here is the key observation.

\begin{proposition}
We have an isomorphism
$$A_u(n)/J=A_o(n)$$
where $J$ is the ideal generated by the relations $v_{ij}=v_{ij}^*$.
\end{proposition}

\begin{proof}
This is clear from definitions of $A_o(n)$ and $A_u(n)$.
\end{proof}

Let $F$ be the set of words on two letters $\alpha,\beta$. For $a\in F$ we denote by $v^{\otimes a}$ the corresponding tensor product of $v=v^{\otimes\alpha}$ and $\bar{v}=v^{\otimes \beta}$.

We denote as usual by $u$ the fundamental corepresentation of $A_o(n)$. Since morphisms increase Hom spaces, we have inclusions
$$Hom(v^a,v^b)\subset Hom(u^{\otimes l(a)},u^{\otimes l(b)})$$
where $l$ is the length of words. These can be combined with equalities in theorem 3.1. We get in this way inclusions
$$Hom(v^a,v^b)\subset TL(l(a),l(b)).$$

\begin{definition}
For $a,b\in F$ we consider the subset
$$D(a,b)\subset D(l(a),l(b))$$
consisting of diagrams $p$ such that when putting $a,b$ on points of $p$, each string joins an $\alpha$ letter to a $\beta$ letter.
\end{definition}

In other words, the set $D(a,b)$ can be described as
$$D(a,b)=\left\{ \begin{matrix}\cdot\,\cdot\,\cdot & \leftarrow &
    \mbox{word }a\cr W & \leftarrow &
    \mbox{ uncolorable strings}\cr \cdot\,\cdot\,\cdot\,\cdot\,\cdot& \leftarrow &
    \mbox{word }b\end{matrix}\right\}$$
where capital letters denote diagrams formed by non-crossing strings, which cannot be colored $\alpha$ or $\beta$, as to match colors of endpoints.

Consider also the subspace
$$TL(a,b)\subset TL(l(a),l(b))$$
generated by diagrams in $D(a,b)$.

\begin{theorem}
We have an equality of vector spaces
$$Hom(v^{\otimes a},v^{\otimes b})=TL(a,b)$$
where $TL(a,b)$ is identified with its image in $Hom(V^{\otimes l(a)},V^{\otimes l(b)})$.
\end{theorem}

\begin{proof}
We follow the proof of theorem 3.1, with notations from there. The starting remark is that for a unitary matrix $v$, the fact that $v^t$ is unitary is equivalent to the fact that $\xi$ is fixed by both $v\otimes\bar{v}$ and $\bar{v}\otimes v$. In other words, we have the following two conditions:
$$E\in Hom(1,v\otimes \bar{v})$$
$$E\in Hom(1,\bar{v}\otimes v)$$

Now since $E$ is the semicircle in $D(0,2)$, these conditions are
$$\cap_1\in Hom(1,v^{\otimes\alpha\beta})$$
$$\cap_2\in Hom(1,v^{\otimes\beta\alpha})$$
where $\cap_1$ is the semicircle having endpoints $\alpha,\beta$, and $\cap_2$ is the semicircle having endpoints $\beta,\alpha$. As in proof of theorem 3.1, this gives
$$<\cap_1,\cap_2>=\{ Hom(v^{\otimes a},v^{\otimes b})\mid a,b\}$$
where tensor categories have this time $F$ as monoid of objects. On the other hand, pictures show that we have the equality
$$<\cap_1,\cap_2>=\{ TL(a,b)\mid a,b\}$$
which together with the above equality gives the result.
\end{proof}

Observe that what changed with respect to proof of theorem 3.1 is the fact that the Temperley-Lieb algebra is replaced with a kind of free version of it. The whole combinatorics is worked out in detail in \cite{ba1}.

We get another proof of a main result in \cite{ba1}, a bit more enlightening than the original one. For two other proofs, probably even more enlightening, but relying on quite technical notions, see \cite{ba12} and \cite{bi}.

\begin{theorem}
We have an embedding of reduced Hopf algebras
$$A_u(n)_{red}\subset {\mathbb C}^*({\mathbb Z})*_{red}A_o(n)_{red}$$
given by $v=zu$, where $z$ is the generator of ${\mathbb Z}$.
\end{theorem}

\begin{proof}
Since $u$ and $u^t$ are unitaries, so are the matrices
$$w=zu$$
$$w^t=zu^t$$
so we get a morphism from left to right:
$$f:A_u(n)\to {\mathbb C}^*({\mathbb Z})*A_o(n)$$

As for any morphism, $f$ increases spaces of fixed points:
$$Hom(1,v^{\otimes a})\subset Hom(1,w^{\otimes a})$$

By standard results in \cite{wo2}, generalising Peter-Weyl theory, $f$ is an isomorphism at level of reduced algebras if and only if all inclusions are equalities. See e.g. \cite{ba1}. Now all fixed point spaces being finite dimensional, this is the same as asking for equalities of dimensions:
$$dim(Hom(1,v^{\otimes a}))= dim(Hom(1,w^{\otimes a}))$$

In terms of characters, we have to prove the formula
$$\int\chi(v)^a=\int\chi(w)^a$$
where exponentials $x^a$ are obtained as corresponding products of terms $x^\alpha=x$ and $x^\beta=x^*$. The term on the right is the $a$-th moment of 
$$\chi(w)=\chi(zu)=z\chi(u)$$
which by Voiculescu's polar decomposition result in \cite{vo} is a circular variable. As for the term on the left, this is given by
$$\int\chi(v)^a=dim(Hom(1,v^{\otimes a}))=\# D(a)$$
which by results of Speicher (\cite{s1}) and Nica-Speicher (\cite{ns}) is also the $a$-th moment of the circular variable.
\end{proof}

\begin{definition}
The $u_{sn}$ variable is given by
$$u_{sn}=v_{11}+\ldots +v_{ss}$$
where $v$ is the fundamental corepresentation of $A_u(n)$.
\end{definition}

This notation looks a bit confusing, because $u_{ij}$ was so far reserved for the fundamental corepresentation of $A_o(n)$. However, this corepresentation will no longer appear, and there is no confusion.

The properties of $u_{sn}$ can be deduced from corresponding properties of $o_{sn}$ by using standard free probability tools.

\begin{theorem}
The $u_{sn}$ variable has the following properties.

(1) We have $u_{sn}=zo_{sn}$, where $z$ is a Haar-unitary free from $o_{sn}$.

(2) The variable $u_{nn}$ is circular.

(3) The variable $(n/s)^{1/2}u_{sn}$ with $n\to\infty$ is circular.
\end{theorem}

\begin{proof}
The first assertion follows from theorem 9.2. The other ones follow from (1) and from corollaries 5.1 and 5.2, by using Voiculescu's result on the polar decomposition of circular variables (\cite{vo}). 
\end{proof}

\begin{theorem}
The elements $(n^{1/2}v_{ij})_{i,j=1,\ldots ,s}$ of $A_u(n)$ with $n\geq s$ become asymptotically free and circular as $n\to\infty$.
\end{theorem}

\begin{proof}
This follows along the same lines as theorem 6.1.
\end{proof}

Finally, it is possible to derive from theorem 9.1 a general integration formula for $A_u(n)$, in the same way as theorem 4.1 is derived from theorem 3.1. For this purpose, we first extend definition 4.2.

\begin{definition}
For $a\in F$, the Gram and Weingarten matrices are
$$G_{an}(p,q)=n^{l(p,q)}$$
$$W_{an}=G_{an}^{-1}$$
where both indices $p,q$ are diagrams in $D(a)$.
\end{definition}

It is convenient at this point to remove the tensor sign in our notations 
$v=v^{\otimes \alpha}$ and $\bar{v}=v^{\otimes \beta}$. 
That is, we use the following notations:
$$v=v^{\alpha}$$
$$\bar{v}=v^{\beta}$$

As in case of $A_o(n)$, we get that integrals are either $0$, or equal to certain sums of entries of the Weingarten matrix.

\begin{theorem}
The Haar functional of $A_u(n)$ is given by
$$\int v_{i_{1}j_{1}}^{a_1}\ldots v_{i_{2k}j_{2k}}^{a_{2k}}=\sum_{pq}
\delta_{pi}\delta_{qj}W_{an}(p,q)$$
$$\int v_{i_{1}j_{1}}^{a_1}\ldots v_{i_{l}j_{l}}^{a_l}=0$$
where $a=a_1a_2\ldots$ is a word in $F$, which in the first formula contains as many $\alpha$ as many $\beta$, and in the second formula, doesn't.
\end{theorem}

\begin{proof}
This proof is done along the same lines as the proof of Theorem 4.1.
%

\end{proof}

Theorem 9.5 has its own interest; however, it is not really needed for study of $u_{sn}$, where the procedure to follow is explained in theorem 9.3 and its proof: find results about $o_{sn}$, then make a free convolution by a Haar-unitary. This kind of convolution operation is standard in free probability, see for instance Nica and Speicher (\cite{ns}).

\end{document}